\documentclass{m2an}
\usepackage{amsmath}
\usepackage{amssymb}
\usepackage{amsthm}
\usepackage{enumerate}
\usepackage{graphicx}
\usepackage{epstopdf}
\usepackage{subfigure}
\usepackage[top=1in, bottom=1in, left=1in, right=1in]{geometry}
\usepackage{tabularx}
\usepackage{changepage}
\usepackage[square,numbers,sort&compress]{natbib}
\usepackage{enumitem}

\newtheorem{theorem}{Theorem}[section]
\newtheorem{lemma}[theorem]{Lemma}

\newtheorem{corollary}[theorem]{Corollary}

\graphicspath{{./Figures/}}

\begin{document}

\title{Supercloseness of Orthogonal Projections onto Nearby Finite Element Spaces}

\author{Evan S. Gawlik}\address{Computational and Mathematical Engineering, Stanford University, Stanford, CA, USA; \email{egawlik@stanford.edu}}
\author{Adrian J. Lew}\secondaddress{Mechanical Engineering, Stanford University, Stanford, CA, USA; \email{lewa@stanford.edu}}\sameaddress{1}


\date{  }

\begin{abstract}
We derive upper bounds on the difference between the orthogonal projections of a smooth function $u$ onto two finite element spaces that are nearby, in the sense that the support of every shape function belonging to one but not both of the spaces is contained in a common region whose measure tends to zero under mesh refinement.  The bounds apply, in particular, to the setting in which the two finite element spaces consist of continuous functions that are elementwise polynomials over shape-regular, quasi-uniform meshes that coincide except on a region of measure $O(h^\gamma)$, where $\gamma$ is a nonnegative scalar and $h$ is the mesh spacing.  The projector may be, for example, the orthogonal projector with respect to the $L^2$- or $H^1$-inner product.  In these and other circumstances, the bounds are superconvergent under a few mild regularity assumptions.  That is, under mesh refinement, the two projections differ in norm by an amount that decays to zero at a faster rate than the amounts by which each projection differs from $u$.  We present numerical examples to illustrate these superconvergent estimates and verify the necessity of the regularity assumptions on $u$.
\end{abstract}



\maketitle

\makeatletter
\renewcommand\paragraph[1]{{\bf #1} }
\makeatother

\section{Introduction} \label{sec:intro}

One of the hallmarks of the finite element method is its geometric flexibility: it permits the construction of numerical approximations to solutions of partial differential equations using meshes that are designed according to the practitioner's discretion.  When two meshes are used to solve the same problem, the norm of the difference between the corresponding numerical solutions is, of course, no larger than the sum of the norms of the differences between each numerical solution and the exact solution.  This paper addresses the question of whether or not a sharper estimate holds in the event that the two meshes coincide over a large fraction of the domain.

Beyond its inherent mathematical appeal, the question raised above has important consequences in the study of numerical solutions to time-dependent PDEs on meshes that change abruptly in time.  Notable examples are remeshing during finite element simulations of problems with moving boundaries, and adaptive refinement during finite element simulations of problems on fixed (or moving) domains.  The relevance of the aforementioned question in this setting is elucidated in~\cite{Gawlik2012b}, where it is shown that if a parabolic PDE is discretized in space with finite elements and the solution is transferred finitely many times between meshes using a suitable projector, then it is possible to derive an upper bound on the error in the numerical solution at a fixed time $T>0$ that involves the norms of the jumps in $r_h u(t)$ across the remeshing times, where $r_h u(t)$ denotes an elliptic projection of the exact solution $u(t)$ onto the current finite element space.  These jumps are precisely the differences between the finite element solutions of an elliptic PDE on two different meshes.

\paragraph{Intuition.}%
It is perhaps not surprising that two finite element solutions associated with nearly identical meshes should differ by an amount that is small relative to their individual differences with the exact solution, under suitable conditions on the finite element spaces and the PDE under consideration.  To develop some intuition, it is instructive to first consider the similarity between the \emph{interpolants} of a smooth function $u$ onto two finite element spaces associated with nearby meshes.  

To this end, consider two families of shape-regular, quasi-uniform meshes $\{\mathcal{T}_h\}_{h \le h_0}$ and $\{\mathcal{T}_h^+\}_{h \le h_0}$ of an open, bounded, Lipschitz domain $\Omega \subset \mathbb{R}^d$, $d \ge 1$.  Assume that the two families are parametrized by a scalar $h$ that equals the maximum diameter of an element among all elements of $\mathcal{T}_h$ and $\mathcal{T}_h^+$ for every $h \le h_0$, where $h_0$ is a positive scalar.  Let $\mathcal{V}_h$ and $\mathcal{V}_h^+$ be finite element spaces consisting of, for definiteness, continuous functions that are elementwise polynomials of degree at most $r-1$ over $\mathcal{T}_h$ and $\mathcal{T}_h^+$, respectively, where $r>1$ is an integer.

For $s \ge 0$ and $p \in [1,\infty]$, we denote by $W^{s,p}(\Omega)$ the Sobolev space of differentiability $s$ and integrability $p$, equipped with the norm $\|\cdot\|_{s,p}$ and semi-norm $|\cdot|_{s,p}$.  We sometimes write $\|\cdot\|_{s,p,\Omega}$ and $|\cdot|_{s,p,\Omega}$ to emphasize the domain under consideration.  We denote $H^s(\Omega) = W^{s,2}(\Omega)$ for every $s \ge 1$ and $L^p(\Omega) = W^{0,p}(\Omega)$ for every $p \in [1,\infty]$.

For finite element spaces of the aforementioned type, the nodal interpolants $i_h u \in \mathcal{V}_h$ and $i_h^+ u \in \mathcal{V}_h^+$ of a function $u \in W^{r,\eta}(\Omega) \cap C^0(\overline{\Omega})$ onto $\mathcal{V}_h$ and $\mathcal{V}_h^+$, respectively, satisfy the standard interpolation estimate
\begin{equation} \label{standard_interp}
\|i_h^+ u - u\|_{s,\eta} + \|i_h u - u\|_{s,\eta} \le C h^{r-s} |u|_{r,\eta}
\end{equation}
for any $s \in \{0,1\}$, any $\eta \in [2,\infty]$, and every $h \le h_0$~\cite{Ern2004}.  Here and throughout this paper, the letter $C$ denotes a constant that is not necessarily the same at each occurrence and is independent of $h$.  

Using the triangle inequality and~(\ref{standard_interp}) with $\eta=2$ gives an immediate upper bound on the $L^2$- and $H^1$-norms of the difference between $i_h^+ u$ and $i_h u$.  Namely,
\begin{equation} \label{naive_estimate}
\|i_h^+ u - i_h u \|_{s,2} \le C h^{r-s} |u|_{r,2}
\end{equation}
for any $s \in \{0,1\}$ and every $h \le h_0$.

Suppose, however, that $\mathcal{T}_h$ and $\mathcal{T}_h^+$ are nearby in the following sense: the two meshes coincide except on a region of measure $O(h^\gamma)$ for some scalar $\gamma \ge 0$. In this scenario, $i_h u$ and $i_h^+ u$ agree everywhere except in the region over which the meshes differ.  Hence, by an application of Holder's inequality (cf. Lemma~\ref{lemma:supp}), the triangle inequality, and~(\ref{standard_interp}),
\begin{align}
\|i_h^+ u - i_h u\|_{s,2} 
&\le C h^{\gamma(1/2-1/\eta)} \|i_h^+ u - i_h u\|_{s,\eta} \nonumber \\
&\le C h^{\gamma(1/2-1/\eta)} \left( \|i_h^+ u - u\|_{s,\eta} + \|u - i_h u\|_{s,\eta}  \right) \nonumber \\
&\le C h^{r-s+\gamma(1/2-1/\eta)} |u|_{r,\eta} \label{interp_superclose}
\end{align}
for any $s \in \{0,1\}$, any $\eta \in [2,\infty]$, and every $h \le h_0$.

A comparison of~(\ref{interp_superclose}) with the naive estimate~(\ref{naive_estimate}) reveals that $i_h u$ and $i_h^+ u$ are \emph{superclose} in the $L^2$- and $H^1$-norms when the corresponding meshes are nearby.  The primary goal of this paper is to prove an analogous superconvergent estimate when $i_h u$ and $i_h^+ u$ are replaced by the orthogonal projections $r_h u$ and $r_h^+ u$ of $u$ onto $\mathcal{V}_h$ and $\mathcal{V}_h^+$, respectively, with respect to a coercive, continuous bilinear form $a : \mathcal{V} \times \mathcal{V} \rightarrow \mathbb{R}$, where $\mathcal{V} \subseteq H^s(\Omega)$ and $s$ is a nonnegative integer.  As special cases, our results apply to $L^2$-projections (the case $s=0$) and elliptic projections (the case $s=1$) onto piecewise polynomial finite element spaces.  Another applicable case of interest is that in which the bilinear form $a$ is of the form 
\[
a(u,w) = \int_\Omega \nabla u \cdot \nabla w \, dx - \int_\Omega (v \cdot \nabla u) w \, dx + \kappa \int_\Omega u w \, dx
\]
with a constant $\kappa>0$ and a vector field $v : \Omega \rightarrow \mathbb{R}^d$.  This bilinear form appears in the analysis of finite element methods for the diffusion equation on a moving domain~\cite{Gawlik2012b}, with $v$ playing the role of the velocity of a moving mesh and $\kappa$ an auxiliary constant introduced to ensure coercivity.

It is not obvious that superconvergent estimates of the form~(\ref{interp_superclose}) should hold in these settings, since the projections of $u$ onto $\mathcal{V}_h$ and $\mathcal{V}_h^+$ need not agree on the region over which the meshes coincide.  Nevertheless, Corollaries~\ref{corollary1} and~\ref{corollary2} provide such estimates under suitable assumptions on the finite element spaces $\mathcal{V}_h$ and $\mathcal{V}_h^+$ and the bilinear form $a$.  The proof uses the observation that, loosely speaking, $a(r_h^+ u - r_h u, r_h^+ u - r_h u)$ is small if $r_h^+ u - r_h u$ is well-approximated by an element of $\mathcal{V}_h^+ \cap \mathcal{V}_h$, since
\[
a(r_h^+ u - r_h u, w_h) = a(r_h^+ u - u, w_h) + a(u - r_h u, w_h) = 0
\]
for any $w_h \in \mathcal{V}_h^+ \cap \mathcal{V}_h$.  In particular, if $\|r_h^+ u - r_h u - w_h\|_{s,2}$ decays to zero more rapidly as $h \rightarrow 0$ than do $\|r_h^+ u - u\|_{s,2}$ and $\|r_h u - u\|_{s,2}$, then a superconvergent estimate for $\|r_h^+ u - r_h u\|_{s,2}$ follows from the relation
\[
a(r_h^+ u - r_h u, r_h^+ u - r_h u) = a(r_h^+ u - r_h u, r_h^+ u - r_h u - w_h)
\]
together with the coercivity and continuity of $a$.  We in fact prove a more general result that applies to the case in which the projectors $r_h$ and $r_h^+$ are associated not only with different subspaces $\mathcal{V}_h$ and $\mathcal{V}_h^+$, but also with different bilinear forms $a_h$ and $a_h^+$ that may depend on $h$.  

\paragraph{Organization.}%
This paper is organized as follows.  In Section~\ref{sec:summary}, we summarize our main results.  We begin with an abstract estimate (Theorem~\ref{thm:abstract_estimate}) for the $H^s$-norm of $r_h^+ u - r_h u$.  We then apply Theorem~\ref{thm:abstract_estimate} to the setting of finite element spaces with nontrivial intersection in Theorem~\ref{thm:Hsestimate}.  Under some additional assumptions on the finite element spaces, the bilinear forms, and the regularity of $u$, we deduce in Corollary~\ref{corollary1} a superconvergent estimate for $\|r_h^+ u - r_h u\|_{s,2}$ that parallels~(\ref{interp_superclose}).  Next, we specialize to the case in which $s=1$ and $a_h$ and $a_h^+$ are bilinear forms associated with elliptic operators that possess smoothing properties.  We use a duality argument to prove a superconvergent estimate (Theorem~\ref{L2estimate} and Corollary~\ref{corollary2}) for the $L^2$-norm of $r_h^+ u - r_h u$ that is up to one order higher than the corresponding estimate in the $H^1$-norm given by Corollary~\ref{corollary1}.

In Section~\ref{sec:proofs}, we present proofs of the preceding results and provide a few remarks along the way.

In Section~\ref{sec:regularity}, we demonstrate the necessity of the regularity assumptions on $u$ that are imposed in the theorems by exhibiting an example of a pair of projectors $r_h$ and $r_h^+$ and a function $u$ whose insufficient regularity leads to a reduction in the rates of convergence of $\|r_h^+u - r_h u\|_{1,2}$ and $\|r_h^+u - r_h u\|_{0,2}$. 

Finally, we give numerical examples to illustrate our positive theoretical results in Section~\ref{sec:numerical}.

\paragraph{Related work.}%
The results presented in this paper bear resemblance to the well-studied phenomenon of superconvergence in finite element theory, where the functions under comparison are typically the solution to a PDE and the numerical solution to a finite element discretization of the same problem.
The phenomenon often manifests itself as an exceptional rate of convergence of the finite element solution to the exact solution at isolated points in the domain, as in~\cite{Barlow1976,Kvrivzek1987,Babuska1996,Goodsell1994,Schatz1996,Wahlbin1995}.  Related results involve exceptional rates of convergence of the finite element solution to a discrete representative of the exact solution, such as its interpolant~\cite{Oganesyan1969,Li2004,Huang2008,Liu2012,Bank2003,Andreev2005,Brandts2003}.  Finally, post-processing techniques can lead to modifications of a finite element solution that converge more rapidly to the exact solution than the unprocessed finite element solution~\cite{Zienkiewicz1992,Babuska1984,Bank2003,Kvrivzek1984,Kvrivzek1987,Goodsell1989,Cockburn2003}.  To our knowledge, however, little attention has been paid to the supercloseness of finite element solutions associated with differing meshes.

\section{Statement of Results} \label{sec:summary}

\paragraph{Notation.}%
Fixing a nonnegative integer $s$ and an open, bounded, Lipschitz domain $\Omega \subset \mathbb{R}^d$, let $\mathcal{V}$ be a closed subspace of $H^s(\Omega)$.  Let $a_h : \mathcal{V} \times \mathcal{V} \rightarrow \mathbb{R}$ and $a_h^+ : \mathcal{V} \times \mathcal{V} \rightarrow \mathbb{R}$ be bilinear forms that may depend on a parameter $h \le h_0$, where $h_0$ is a positive scalar.  We assume that $a_h$ and $a_h^+$ are continuous and coercive uniformly in $h$.  In other words, for every $h \le h_0$ and every $u,w \in \mathcal{V}$, the inequalities
\begin{align*}
a_h(u,u) &\ge \alpha\|u\|_{s,2}^2, \\
a_h(u,w) &\le M\|u\|_{s,2} \|w\|_{s,2}
\end{align*}
hold with constants $\alpha$ and $M$ independent of $h$, and similarly for $a_h^+$ (with the same constants $\alpha$ and $M$).  

Let $\{\mathcal{V}_h\}_{0<h\le h_0}$ and $\{\mathcal{V}_h^+\}_{0<h\le h_0}$ be two families of finite element subspaces of $\mathcal{V}$.  It is a consequence of the Lax-Milgram theorem that the maps $r_h : \mathcal{V} \rightarrow \mathcal{V}_h$ and $r_h^+ : \mathcal{V} \rightarrow \mathcal{V}_h^+$ defined by the relations
\begin{equation*}
a_h(r_h u - u, w_h) = 0 \;\;\; \forall w_h \in \mathcal{V}_h
\end{equation*}
and
\begin{equation*}
a_h^+(r_h^+ u - u, w_h^+) = 0 \;\;\; \forall w_h^+ \in \mathcal{V}_h^+,
\end{equation*}
respectively, are well-defined linear projectors.  

For intuition, it is useful to think of $\mathcal{V}_h$ and $\mathcal{V}_h^+$ as finite element spaces associated with a pair of meshes $\mathcal{T}_h$ and $\mathcal{T}_h^+$ of $\Omega$, with the parameter $h$ denoting the maximum diameter of an element among all elements of $\mathcal{T}_h$ and $\mathcal{T}_h^+$.  This level of concreteness, however, is not needed for a presentation of the results that follow.

\paragraph{Abstract estimate.}%
Our first result is an abstract estimate for the $H^s$-norm of $r_h^+ u - r_h u$.  It provides an alternative to the obvious upper bound
\[
\|r_h^+ u - r_h u\|_{s,2} \le \|r_h^+ u - u\|_{s,2} + \|u - r_h u\|_{s,2}
\]
that one obtains from the triangle inequality.  Its utility will be made apparent shortly.

\begin{theorem} \label{thm:abstract_estimate}
Let $a_h^+$ and $a_h$ be uniformly coercive and continuous bilinear forms on $\mathcal{V} \times \mathcal{V}$.  Then for every $u \in \mathcal{V}$ and every $h \le h_0$,
\begin{equation} \label{abstract_estimate}
\begin{split}
\| r_h^+ u - r_h u \|_{s,2} \le \inf_{\substack{e_h \in \mathcal{V}_h \\ e_h^+ \in \mathcal{V}_h^+}} 
\Big[ \frac{M}{\alpha} &\left\| r_h^+ u - r_h u - (e_h + e_h^+) \right\|_{s,2} \\
+ &\frac{1}{\sqrt{\alpha}}  \big( \left|a_h^+(r_h^+u - u, e_h)\right|^{1/2} + \left|a_h(r_h u - u, e_h^+)\right|^{1/2} \\
&+ \left|a_h^+(r_h u - u, e_h + e_h^+) - a_h(r_h u - u, e_h + e_h^+)\right|^{1/2} \big) \Big].
\end{split}
\end{equation}
\end{theorem}

The preceding theorem provides a heuristic for estimating the $H^s$-norm of $r_h^+ u - r_h u$.  Namely, one seeks functions $e_h \in \mathcal{V}_h$ and $e_h^+ \in \mathcal{V}_h^+$ that are nearly (right-) orthogonal to $r_h^+ u - u$ and $r_h u - u$ with respect to $a_h^+(\cdot,\cdot)$ and $a_h(\cdot,\cdot)$, respectively, but whose sum is close to $r_h^+ u - r_h u$.  In general, near orthogonality and closeness to $r_h^+ u - r_h u$ are competing interests.  Exact orthogonality holds for $e_h,e_h^+ \in \mathcal{V}_h^+ \cap \mathcal{V}_h$, whereas $e_h + e_h^+$ can be made equal to $r_h^+ u - r_h u$ by choosing, for instance, $e_h^+ = r_h^+u$ and $e_h = -r_h u$.  If a suitable choice of $e_h$ and $e_h^+$ leads to adequate satisfaction of both interests simultaneously, and if $a_h^+$ is close to $a_h$ (in the sense that the last term in~(\ref{abstract_estimate}) is small), then the prospects of producing a superconvergent bound on $\|r_h^+ u - r_h u\|_{s,2}$ are favorable.

\paragraph{Finite element spaces with nontrivial intersection.}%
We now apply Theorem~\ref{thm:abstract_estimate} to the case in which the finite element spaces $\mathcal{V}_h^+$ and $\mathcal{V}_h$ intersect nontrivially.  The setting that we have in mind is that in which $\mathcal{V}_h$ and $\mathcal{V}_h^+$ consist of continuous functions that are elementwise polynomials over shape-regular, quasi-uniform meshes of $\Omega$ that coincide except on a region of measure $O(h^\gamma)$ for some constant $\gamma \ge 0$.  To allow for more generality, we state the assumptions on $\mathcal{V}_h^+$ and $\mathcal{V}_h$ abstractly, and we refer the reader to Appendix~\ref{sec:appendix_Pk} for a proof of their satisfaction in the aforementioned setting.

In particular, we assume the existence of a constant $\eta \in [2,\infty]$ such that the following properties hold:
\begin{enumerate}[label=(\ref*{thm:Hsestimate}.\roman*),ref=\ref*{thm:Hsestimate}.\roman*]
\item \label{Wseta} For every $h \le h_0$, $\mathcal{V}_h,\mathcal{V}_h^+ \subset W^{s,\eta}(\Omega) \cap \mathcal{V}$.
\item \label{inverse} There exists $C>0$ independent of $h$ such that the inverse estimate
\[
\|w_h\|_{m,\eta} \le Ch^{-m}\|w_h\|_{0,\eta}
\]
holds for every $m=0,1,\dots,s$, every $w_h \in \mathcal{V}_h^+ \cap \mathcal{V}_h$, and every $h \le h_0$.
\item \label{assumption2b} There exist constants $\gamma \ge 0$ and $C>0$ independent of $h$ and a map $\pi_h : \mathcal{V}_h^+ + \mathcal{V}_h \rightarrow \mathcal{V}_h^+ \cap \mathcal{V}_h$ such that
\begin{equation*} 
\|\pi_h w_h\|_{0,\eta} \le C \|w_h\|_{0,\eta}
\end{equation*}
and
\begin{equation*}
|\mathrm{supp}(\pi_h w_h - w_h)| \le C h^\gamma
\end{equation*}
for every $w_h \in \mathcal{V}_h^+ + \mathcal{V}_h$ and every $h \le h_0$.
\end{enumerate}

In the context of finite element spaces consisting of continuous functions that are elementwise polynomials over shape-regular, quasi-uniform meshes of $\Omega$, a befitting choice for $\pi_h$ in~(\ref{assumption2b}) is the nodal interpolant onto $\mathcal{V}_h^+ \cap \mathcal{V}_h$; see Appendix~\ref{sec:appendix_Pk}.  In that setting, the constant $\gamma$ appearing in~(\ref{assumption2b}) may take on any real value between $0$ and $d$, unless the two meshes coincide entirely (in which case $\gamma$ may be chosen arbitrarily large).  To realize a pair of meshes $\mathcal{T}_h$ and $\mathcal{T}_h^+$ fulfilling~(\ref{assumption2b}) with $\gamma \in [0,d]$, one may, for instance, consider a shape-regular, quasi-uniform mesh $\mathcal{T}_h$ of $\Omega$ and perturb the positions of $O(h^{-d+\gamma})$ of its nodes by a sufficiently small amount to define $\mathcal{T}_h^+$.

The following theorem results from applying Theorem~\ref{thm:abstract_estimate} to the setting delineated in conditions~(\ref{Wseta}-\ref{assumption2b}), with the choice $e_h=\pi_h(r_h^+ u - r_h u)$ and $e_h^+=0$ in~(\ref{abstract_estimate}).

\begin{theorem} \label{thm:Hsestimate}
Suppose the conditions of Theorem~\ref{thm:abstract_estimate} hold and the finite element spaces $\mathcal{V}_h^+$ and $\mathcal{V}_h$ satisfy conditions~(\ref{Wseta}-\ref{assumption2b}).  Suppose further that there exist constants $C_1>0$, $\delta \ge 0$, $1 \le q \le \eta$, and $\mu,\nu \in \{0,1,\dots,s\}$ independent of $h$ such that
\begin{equation} \label{assumption1}
|a_h^+(v,w) - a_h(v,w)| \le C_1 h^\delta \|v\|_{\mu,\eta} \|w\|_{\nu,q}
\end{equation}
for every $v,w \in W^{s,\eta}(\Omega) \cap \mathcal{V}$ and every $h \le h_0$.  
Then there exists $C>0$ independent of $h$ such that for any $h \le h_0$ and any $u \in W^{s,\eta}(\Omega) \cap \mathcal{V}$,
\[
\begin{split}
\|r_h^+ u - r_h u\|_{s,2} \le C h^{\sigma-s} \Big[ &h^s\|r_h^+ u - u\|_{s,\eta} + h^s\|r_h u - u\|_{s,\eta} + \|r_h^+ u - u\|_{0,\eta} + \|r_h u - u\|_{0,\eta} \\
&+ \left(h^\mu \|r_h u - u\|_{\mu,\eta}\right)^{1/2} \left(\|r_h^+ u - u\|_{0,\eta} + \|r_h u - u\|_{0,\eta} \right)^{1/2} \Big]
\end{split}
\]
with
\begin{equation} \label{sigma}
\sigma = \min\left\{\gamma\left(\frac{1}{2}-\frac{1}{\eta}\right), \frac{\delta+2s-\mu-\nu}{2} \right\}.
\end{equation}
\end{theorem}

The meaning of Theorem~\ref{thm:Hsestimate} is clearest when the quantities $h^m \|r_h u - u\|_{m,p}$ and $h^m \|r_h^+ u - u\|_{m,p}$, $m = 0,1,\dots,s$, $p=2,\eta$, all decay at the same rate with respect to $h$ as $h \rightarrow 0$.  In such a setting, the theorem states that $\|r_h^+ u - r_h u\|_{s,2}$ tends to zero faster than $\|r_h u - u\|_{s,2} + \|r_h^+ u - u\|_{s,2}$ by a factor $O(h^\sigma)$, where the order of superconvergence $\sigma$ depends primarily upon two features: (1) the extent to which the finite element spaces $\mathcal{V}_h$ and $\mathcal{V}_h^+$ coincide, as measured by the constant $\gamma$ in~(\ref{assumption2b}), and (2) the difference between the bilinear forms $a_h$ and $a_h^+$, as measured by the constants $\delta$, $\mu$, and $\nu$ in~(\ref{assumption1}).  The regularity of $u$ also plays a role in the estimate via the constant $\eta$, which is in the best case equal to $\infty$.

To be more concrete, let us point out that in many contexts (which we detail in Appendix~\ref{sec:appendix_Linf}), the quantities $r_h u - u$ and $r_h^+ u - u$ satisfy estimates of the form
\begin{align}
\|r_h u - u\|_{0,\eta} + \|r_h^+ u - u\|_{0,\eta} &\le C \ell(h) h^r |u|_{r,\eta} \label{ep0eta}, \\
\|r_h u - u\|_{m,\eta} + \|r_h^+ u - u\|_{m,\eta} &\le C h^{r-m} |u|_{r,\eta}, \;\;\; m=1,2,\dots,s, \label{epmeta}
\end{align}
for every $u \in W^{r,\eta}(\Omega) \cap \mathcal{V}$ and every $h \le h_0$, where $r>s$ is an integer and $\ell(h)$ is either identically unity or equal to $\log(h^{-1})$.  Note that~(\ref{epmeta}) is vacuous when $s=0$.  When such estimates hold, the following corollary to Theorem~\ref{thm:Hsestimate} is immediate.
\begin{corollary} \label{corollary1}
Suppose that the conditions of Theorem~\ref{thm:Hsestimate} are satisfied and that both $r_h$ and $r_h^+$ satisfy estimates of the form~(\ref{ep0eta}-\ref{epmeta}) for an integer $r>s$.  Then there exists $C>0$ independent of $h$ such that
\begin{equation*}
\|r_h^+ u - r_h u\|_{s,2} \le C \ell(h) h^{r-s+\sigma} |u|_{r,\eta}
\end{equation*}
for every $u \in W^{r,\eta}(\Omega) \cap \mathcal{V}$ and every $h \le h_0$, with $\sigma$ given by~(\ref{sigma}).

In particular, if $a_h=a_h^+$, then
\begin{equation*}
\|r_h^+ u - r_h u\|_{s,2} \le C \ell(h) h^{r-s+\gamma(1/2-1/\eta)} |u|_{r,\eta}
\end{equation*}
for every $u \in W^{r,\eta}(\Omega) \cap \mathcal{V}$ and every $h \le h_0$.
\end{corollary}

Note that to deduce the preceding corollary, the case $a_h=a_h^+$ is handled by taking $\delta = \infty$ and choosing any admissible $\mu$, $\nu$ and $q$ in~(\ref{assumption1}).

\paragraph{$L^2$ estimates for elliptic projections.}%
Finally, we restrict our attention to the case $s=1$ with $\mathcal{V}=H^1_0(\Omega)$, so that $a_h$ and $a_h^+$ are coercive, continuous bilinear forms on $H^1_0(\Omega) \times H^1_0(\Omega)$, uniformly in $h$.  Here, $H^1_0(\Omega)$ denotes the space of functions in $H^1(\Omega)$ with vanishing trace on $\partial \Omega$. Our aim is to provide an estimate for the $L^2$-norm of $r_h^+ u - r_h u$ that parallels the estimate in the $H^1$-norm provided by Corollary~\ref{corollary1} but is of a higher order by up to one power of $h$.  

In addition to the assumptions stated in Theorem~\ref{thm:Hsestimate}, we make the following assumptions on the bilinear forms $a_h$ and $a_h^+$.

\begin{enumerate}[label=(\ref*{L2estimate}.\roman*),ref=\ref*{L2estimate}.\roman*]
\item \label{cond:smoothing} The bilinear forms $a_h$ and $a_h^+$ are associated with elliptic operators whose adjoints possess \emph{smoothing properties} (cf.~\cite[Definition 3.14]{Ern2004}), uniformly in $h$.  Precisely, let $f \in L^2(\Omega)$ and consider the following problem: Find $w \in \mathcal{V}$ such that
\begin{equation} \label{adjoint_problem}
a_h(y,w) = (f,y) \;\;\; \forall y \in \mathcal{V}, 
\end{equation}
where $(f,y) := \int_\Omega fy$.  Then $a_h$ is said to have smoothing properties (uniformly in $h$) if there exists a constant $C>0$ independent of $h$ such that for every $f \in L^2(\Omega)$ and every $h \le h_0$, there exists a unique solution $w$ to~(\ref{adjoint_problem}) satisfying the elliptic regularity estimate
\begin{equation*} 
\|w\|_{2,2} \le C\|f\|_{0,2}.
\end{equation*}
\item \label{cond:subdomains} There exists $C>0$ such that for any $h \le h_0$, any subdomain $R \subseteq \Omega$, and any $v,w \in \mathcal{V}$ with $\mathrm{supp}(w) \subseteq R$,
\begin{equation*} 
|a_h(v,w)| \le C \|v\|_{1,2,R} \|w\|_{1,2,R},
\end{equation*}
where the constant $C$ is independent of $h$ and $R$, and similarly for $a_h^+$.
\item \label{qrestriction} The constant $q$ appearing in the bound~(\ref{assumption1}) satisfies the additional restriction 
\begin{equation*}
\begin{cases}
q < \infty &\mbox{ if } d=4-2\nu, \\
q \le \frac{2d}{d-4+2\nu} &\mbox{ if } d > 4-2\nu.
\end{cases}
\end{equation*}
\end{enumerate}
Condition~(\ref{qrestriction}) guarantees the validity of the Sobolev emdedding $H^2(\Omega) \subset W^{\nu,q}(\Omega)$.  Note that it places no additional restriction on $q$ if $d < 4-2\nu$.

Furthermore, we assume the existence of interpolation operators $i_h : \bar{\mathcal{V}} \rightarrow \mathcal{V}_h$ and $i_h^+ : \bar{\mathcal{V}} \rightarrow \mathcal{V}_h^+$ defined on a space $H^2(\Omega) \cap \mathcal{V} \subseteq \bar{\mathcal{V}} \subseteq \mathcal{V}$ that satisfy the following properties.
\begin{enumerate}[label=(\ref*{L2estimate}.\roman*),ref=\ref*{L2estimate}.\roman*,resume]
\item \label{cond:interp_stability} There exists $C>0$ independent of $h$ such that 
\[
\|i_h w\|_{\nu,q} + \|i_h^+ w\|_{\nu,q} \le C\|w\|_{\nu,q}
\]
for every $w \in H^2(\Omega) \cap \mathcal{V}$ and every $h \le h_0$.
\item \label{cond:interp} There exists $C>0$ independent of $h$ such that
\[
\|i_h w - w\|_{1,2} + \|i_h^+ w - w\|_{1,2} \le C h |w|_{2,2}
\]
for every $w \in H^2(\Omega) \cap \mathcal{V}$ and every $h \le h_0$.
\item \label{cond:interp_coincide} For every $w \in H^2(\Omega) \cap \mathcal{V}$ and every $h \le h_0$, 
\[
\mathrm{supp}(i_h^+ w - i_h w) \subseteq \mathcal{R}_h,
\]
where
\[
\mathcal{R}_h := \bigcup_{w_h \in \mathcal{V}_h + \mathcal{V}_h^+} \mathrm{supp}(w_h - \pi_h w_h)
\]
and $\pi_h$ is the map introduced in~(\ref{assumption2b}).
\end{enumerate}

Our estimate for the $L^2$-norm of $r_h^+ u - r_h u$, whose proof employs a duality argument, is as follows.

\begin{theorem} \label{L2estimate}
Suppose the conditions of Theorem~\ref{thm:Hsestimate} hold with $s=1$. Assume further that conditions~(\ref{cond:smoothing}-\ref{cond:interp_coincide}) hold.  Then there exists $C>0$ independent of $h$ such that for every $u \in W^{1,\eta}(\Omega) \cap \mathcal{V}$ and every $h \le h_0$, 
\[
\begin{split}
\|r_h^+ u - r_h u\|_{0,2} \le C h^{\sigma'} \, \Big[ &h\|r_h^+ u - u\|_{1,\eta} + h\|r_h u - u\|_{1,\eta} + \|r_h^+ u - u\|_{0,\eta} + \|r_h u - u\|_{0,\eta} \\
&+ \left(h^\mu \|r_h u - u\|_{\mu,\eta}\right)^{1/2} \left(\|r_h^+ u - u\|_{0,\eta} + \|r_h u - u\|_{0,\eta} \right)^{1/2} \\
&+ h^\mu \|r_h u - u\|_{\mu,\eta} \Big],
\end{split}
\]
with
\begin{equation} \label{sigmaprime}
\sigma' = \min\left\{\gamma\left(\frac{1}{2}-\frac{1}{\eta}\right), \frac{\delta+2-\mu-\nu}{2}, \delta - \mu \right\}.
\end{equation}
\end{theorem}

Just as in Theorem~\ref{thm:Hsestimate}, the meaning of Theorem~\ref{L2estimate} is clearest when the quantities $h^m \|r_h u - u\|_{m,p}$ and $h^m \|r_h^+ u - u\|_{m,p}$, $m = 0,1,\dots,s$, $p=2,\eta$, all decay at the same rate with respect to $h$ as $h \rightarrow 0$.  In such a setting, Theorem~\ref{L2estimate} states that $\|r_h^+ u - r_h u\|_{0,2}$ tends to zero faster than $\|r_h u - u\|_{0,2} + \|r_h^+ u - u\|_{0,2}$ by a factor $O(h^{\sigma'})$, where the order of superconvergence $\sigma'$ is given by~(\ref{sigmaprime}).  Note that $\sigma' \le \sigma$, where $\sigma$ is the order of superconvergence of the $H^1$-norm of $r_h^+ u - r_h u$ that was provided in Theorem~\ref{thm:Hsestimate}.

Concretely, when estimates of the form~(\ref{ep0eta}-\ref{epmeta}) hold for $u \in W^{r,\eta}(\Omega) \cap \mathcal{V}$ with an integer $r>1$, we arrive immediately at the following corollary to Theorem~\ref{L2estimate}.
\begin{corollary} \label{corollary2}
Suppose that the conditions of Theorem~\ref{thm:Hsestimate} are satisfied and that both $r_h$ and $r_h^+$ satisfy estimates of the form~(\ref{ep0eta}-\ref{epmeta}) for an integer $r>1$.  Then there exists $C>0$ independent of $h$ such that
\begin{equation*}
\|r_h^+ u - r_h u\|_{0,2} \le C \ell(h) h^{r+\sigma'} |u|_{r,\eta}
\end{equation*}
for every $u \in W^{r,\eta}(\Omega) \cap \mathcal{V}$ and every $h \le h_0$, with $\sigma'$ given by~(\ref{sigmaprime}).  

In particular, if $a_h = a_h^+$, then
\begin{equation*}
\|r_h^+ u - r_h u\|_{0,2} \le C \ell(h) h^{r+\gamma(1/2-1/\eta)} |u|_{r,\eta}
\end{equation*}
for every $u \in W^{r,\eta}(\Omega) \cap \mathcal{V}$ and every $h \le h_0$.
\end{corollary}

Note that to deduce the preceding corollary, the case $a_h=a_h^+$ is again handled by taking $\delta = \infty$ and choosing any admissible $\mu$, $\nu$ and $q$ in~(\ref{assumption1}).

\section{Proofs} \label{sec:proofs}

This section presents proofs of Theorems~\ref{thm:abstract_estimate},~\ref{thm:Hsestimate}, and~\ref{L2estimate}.

\proof[Proof of Theorem~\ref{thm:abstract_estimate}] Let $e_h \in \mathcal{V}_h$ and $e_h^+ \in \mathcal{V}_h^+$, and write
\[
\begin{split}
a_h^+(r_h^+ u - r_h u, r_h^+ u - r_h u) =&\; a_h^+\left(r_h^+ u - r_h u, r_h^+ u - r_h u - (e_h + e_h^+)\right) \\
&+ a_h^+(r_h^+ u - r_h u, e_h + e_h^+).
\end{split}
\]
The uniform coercivity and continuity of $a_h^+$ imply
\[
\| r_h^+ u - r_h u \|_{s,2}^2 \le \frac{1}{\alpha} \left( M \| r_h^+ u - r_h u \|_{s,2} \|r_h^+ u - r_h u - (e_h + e_h^+)\|_{s,2} + |a_h^+(r_h^+ u - r_h u, e_h + e_h^+)| \right).
\]
Using the fact that for real numbers $x,a,b \ge 0$, 
\[
x^2 \le a x + b \implies x \le a + \sqrt{b},
\] 
we deduce that
\[
\| r_h^+ u - r_h u \|_{s,2} \le \frac{M}{\alpha} \|r_h^+ u - r_h u - (e_h + e_h^+)\|_{s,2} + \frac{1}{\sqrt{\alpha}} |a_h^+(r_h^+ u - r_h u, e_h + e_h^+)|^{1/2} 
\]
The result will then follow from the identity
\begin{equation} \label{term2expansion}
\begin{split}
a_h^+(r_h^+ u - r_h u, e_h + e_h^+) 
&= a_h^+(r_h^+ u - u, e_h) + a_h(u - r_h u, e_h^+) \\
&\;\;\;+ a_h^+(u - r_h u, e_h + e_h^+) - a_h(u - r_h u, e_h + e_h^+)
\end{split}
\end{equation}
together with the subadditivity of the square root operator.

To prove~(\ref{term2expansion}), use the decomposition $r_h^+ u - r_h u = (r_h^+ u - u) + (u - r_h u)$ to write
\[
a_h^+(r_h^+ u - r_h u, e_h + e_h^+) = a_h^+(r_h^+ u - u, e_h + e_h^+) + a_h^+(u - r_h u, e_h + e_h^+).
\]
Now add and subtract $a_h(u - r_h u, e_h + e_h^+)$ to obtain
\[
\begin{split}
a_h^+(r_h^+ u - r_h u, e_h + e_h^+) 
&= a_h^+(r_h^+ u - u, e_h + e_h^+) + a_h(u - r_h u, e_h + e_h^+) \\
&\;\;\;+ a_h^+(u - r_h u, e_h + e_h^+) - a_h(u - r_h u, e_h + e_h^+).
\end{split}
\] 
Finally, use the definitions of $r_h^+$ and $r_h$ to simplify the first two terms, giving~(\ref{term2expansion}).
\qed

We remark that while the estimate~(\ref{abstract_estimate}) is not symmetric in the ``+'' variables and their unadorned counterparts, it can easily be made symmetric by exchanging the roles of $r_h^+$ and $a_h^+$ with $r_h$ and $a_h$, respectively, and averaging the resulting estimates.  The same holds true for the estimates in Theorems~\ref{thm:Hsestimate} and~\ref{L2estimate}.

We now turn to the proof of Theorem~\ref{thm:Hsestimate}.  We begin with a lemma concerning the relationship between a function's support and its Sobolev norms.

\begin{lemma} \label{lemma:supp}
Let $f \in W^{k,p}(\Omega)$, $k \ge 0$, $p \in [1,\infty]$.  Then for any $1 \le t \le p$,
\[
\|f\|_{k,t} \le |\mathrm{supp}(f)|^{1/t-1/p} \|f\|_{k,p}.
\]
\end{lemma}
\proof Let $\chi : \Omega \rightarrow \{0,1\}$ denote the indicator function for $\mathrm{supp}(f)$.  We have
\begin{align*}
\|f\|_{k,t} 
&= \sum_{|\alpha| \le k} \|\partial^\alpha f\|_{0,t} \\
&= \sum_{|\alpha| \le k} \|\chi \partial^\alpha f\|_{0,t}.
\end{align*}
Now let $\tilde{p} \in [1,\infty]$ be such that $\frac{1}{\tilde{p}}+\frac{1}{p}=\frac{1}{t}$. By Holder's inequality,
\begin{align*}
\|f\|_{k,t} 
&\le \sum_{|\alpha| \le k} \|\chi\|_{0,\tilde{p}} \|\partial^\alpha f\|_{0,p} \\
&= |\mathrm{supp}(f)|^{1/\tilde{p}} \sum_{|\alpha| \le k} \|\partial^\alpha f\|_{0,p} \\
&= |\mathrm{supp}(f)|^{1/t-1/p} \|f\|_{k,p}.
\end{align*}
\qed

The proof of Theorem~\ref{thm:Hsestimate} is as follows.

\proof[Proof of Theorem~\ref{thm:Hsestimate}]
Choose $e_h^+=0$ and $e_h = \pi_h (r_h^+ u - r_h u)$ in~(\ref{abstract_estimate}).  By the stability assumption in~(\ref{assumption2b}), 
\begin{align*}
\|e_h\|_{0,\eta} 
&\le C\|r_h^+ u - r_h u\|_{0,\eta} \\
&\le C \left( \|r_h^+ u - u\|_{0,\eta} + \|u - r_h u\|_{0,\eta} \right).
\end{align*}
Thus, for any $m=0,1,\dots,s$,
\begin{equation} \label{ehmeta}
\| e_h \|_{m,\eta} \le C h^{-m} \left( \|r_h^+ u - u\|_{0,\eta} + \|u - r_h u\|_{0,\eta} \right)
\end{equation} 
by~(\ref{inverse}).  It follows that
\begin{align*}
\|r_h^+ u - r_h u - (e_h + e_h^+) \|_{s,\eta} 
&\le \|r_h^+ u - u\|_{s,\eta} + \|u - r_h u\|_{s,\eta} + \|e_h\|_{s,\eta} + \|e_h^+\|_{s,\eta} \\
&\le C \big( \|r_h^+ u - u\|_{s,\eta} + \|u - r_h u\|_{s,\eta} \\
&\hspace{0.4in}+ h^{-s}\|r_h^+ u - u\|_{0,\eta} + h^{-s} \|u - r_h u\|_{0,\eta} \big).
\end{align*}
Now note that $r_h^+ u - r_h u - (e_h + e_h^+)$ 
has support of measure $O(h^\gamma)$ by~(\ref{assumption2b}).  Consequently, by Lemma~\ref{lemma:supp},
\begin{align}
\|r_h^+ u - r_h u - (e_h + e_h^+) \|_{s,2} 
&\le C h^{\gamma(1/2-1/\eta)} \|r_h^+ u - r_h u - (e_h + e_h^+) \|_{s,\eta} \nonumber \\
&\le C h^{\gamma(1/2-1/\eta)} \big( \|r_h^+ u - u\|_{s,\eta} + \|r_h u - u\|_{s,\eta} \nonumber \\
&\hspace{0.7in} + h^{-s}\|r_h^+ u - u\|_{0,\eta} + h^{-s} \|r_h u - u\|_{0,\eta}\big). \label{term1}
\end{align}
To estimate the remaining terms that appear in~(\ref{abstract_estimate}), note that
\[
a_h^+(r_h^+ u - u, e_h) = 0
\]
since $e_h \in \mathcal{V}_h^+ \cap \mathcal{V}_h \subseteq \mathcal{V}_h^+$, and
\[
a_h(r_h u - u, e_h^+) = 0
\]
since $e_h^+ = 0$.  Finally, using~(\ref{ehmeta}) with $m=\nu$ together with~(\ref{assumption1}) shows that
\begin{align*}
\big|a_h^+(r_h u - u, e_h + &e_h^+) - a_h(r_h u - u, e_h + e_h^+)\big| \\
&\le C h^{\delta} \|r_h u - u\|_{\mu,\eta} \|e_h\|_{\nu,q} \\
&\le C h^{\delta} \|r_h u - u\|_{\mu,\eta} \|e_h\|_{\nu,\eta} \\
&\le C h^{\delta-\nu} \|r_h u - u\|_{\mu,\eta} \left(\|r_h^+ u - u\|_{0,\eta} + \|u - r_h u\|_{0,\eta}\right).
\end{align*}
Taking the square root and adding~(\ref{term1}) proves the claim.
\qed

Note that the preceding proof treats the estimate~(\ref{assumption1}) wastefully when $q<\eta$, in the sense that the ultimate bound on $\|r_h^+ u - r_h u\|_{s,2}$ is unchanged if $q$ is replaced by $\eta$.  The importance of considering scenarios in which $q$ may be chosen less than $\eta$ is made apparent in Theorem~\ref{L2estimate}, where the restriction~(\ref{qrestriction}) is enforced.

With this in mind, we now prove Theorem~\ref{L2estimate}.

\proof[Proof of Theorem~\ref{L2estimate}]
Define $w \in \mathcal{V}$ as the solution to the dual problem
\begin{equation} \label{dual_problem}
a_h^+(y,w) = (r_h^+ u - r_h u, y) \;\;\; \forall y \in \mathcal{V}. 
\end{equation}
Note that $w \in H^2(\Omega) \cap \mathcal{V}$ by~(\ref{cond:smoothing}).

For any $w_h^+ \in \mathcal{V}_h^+$, $w_h \in \mathcal{V}_h$, we have
\begin{align*}
\|r_h^+ u - r_h u\|_{0,2}^2 
&= a_h^+(r_h^+ u - r_h u, w) \\
&= a_h^+(r_h^+ u - r_h u, w - w_h^+) + a_h^+(r_h^+ u - r_h u, w_h^+) \\
&= a_h^+(r_h^+ u - r_h u, w - w_h^+) + a_h^+(u - r_h u, w_h^+) \\
&= a_h^+(r_h^+ u - r_h u, w - w_h^+) + a_h^+(u - r_h u, w_h^+ - w_h) \\
&\;\;\; + a_h^+(u - r_h u, w_h) - a_h(u-r_h u, w_h) \\
&=: T_1 + T_2 + T_3,
\end{align*}
where
\begin{align*}
T_1 &= a_h^+(r_h^+ u - r_h u, w - w_h^+), \\
T_2 &= a_h^+(u - r_h u, w_h^+ - w_h), \\
T_3 &= a_h^+(u - r_h u, w_h) - a_h(u-r_h u, w_h).
\end{align*}

Now choose $w_h^+ = i_h^+ w$ and $w_h = i_h w$ and bound each term separately.  By the continuity of $a_h^+$ and~(\ref{cond:interp}),
\begin{align*}
|T_1| 
&\le C \| r_h^+ u - r_h u\|_{1,2} \| w - w_h^+ \|_{1,2} \\
&\le C h \| r_h^+ u - r_h u\|_{1,2} |w|_{2,2}. \\
\end{align*}
To bound $T_2$, note that $\mathrm{supp}(w_h^+ - w_h) \subseteq \mathcal{R}_h$ has measure $O(h^\gamma)$ by~(\ref{cond:interp_coincide}) and~(\ref{assumption2b}).  Thus,
\begin{align*}
|T_2|
&\le C \|u-r_h u\|_{1,2,\mathcal{R}_h} \|w_h^+ - w_h\|_{1,2,\mathcal{R}_h} \\
&\le C h^{\gamma(1/2-1/\eta)} \|u-r_h u\|_{1,\eta} \left(\|w_h^+ - w\|_{1,2,\mathcal{R}_h} + \|w - w_h\|_{1,2,\mathcal{R}_h} \right)  \\
&\le C h^{\gamma(1/2-1/\eta)+1} \|u-r_h u\|_{1,\eta}  |w|_{2,2}
\end{align*}
by~(\ref{cond:subdomains}), Lemma~\ref{lemma:supp}, and~(\ref{cond:interp}).
For $T_3$, we have by~(\ref{assumption1}) that
\begin{equation*}
|T_3| \le C h^\delta \|u-r_h u\|_{\mu,\eta} \|w_h\|_{\nu,q}.
\end{equation*}
Using~(\ref{cond:interp_stability}) together with the Sobolev embedding $H^2(\Omega) \subset W^{\nu,q}(\Omega)$ ensured by~(\ref{qrestriction}) gives
\begin{equation*}
|T_3| \le C h^\delta \|u-r_h u\|_{\mu,\eta} \|w\|_{2,2}.
\end{equation*}
Combining results and invoking the regularity estimate~(\ref{cond:smoothing}) leads to
\[
\begin{split}
\|r_h^+ u - r_h u\|_{0,2} \le C \Big[ h &\|r_h^+ u - r_h u\|_{1,2} \\&+ h^{\min\{\gamma(1/2-1/\eta),\delta-\mu\}} \left( h \|u-r_h u\|_{1,\eta} + h^\mu \|u-r_h u\|_{\mu,\eta} \right) \Big].
\end{split}
\]
Conclude using Theorem~\ref{thm:Hsestimate}.
\qed


\section{The Need for Regularity} \label{sec:regularity}

When $a_h^+=a_h$ and $\gamma$ is fixed, the estimates of Corollaries~\ref{corollary1} and~\ref{corollary2} are of the highest order in $h$ when $\eta = \infty$, but in this case they demand that $u \in W^{r,\infty}(\Omega) \cap \mathcal{V}$.  If the regularity requirement $u \in W^{r,\infty}(\Omega) \cap \mathcal{V}$ is relaxed, the rates of convergence of $\|r_h^+ u - r_h u\|_{0,2}$ and $\|r_h^+ u - r_h u\|_{1,2}$ as $h \rightarrow 0$ may deteriorate.

Indeed, consider the case in which $\mathcal{V}_h$ is the space of piecewise affine functions on a grid $(0,h,2h,3h,\dots,1)$ of the unit interval in one dimension that vanish at 0 and 1.  Let $\mathcal{V}_h^+$ be the space of piecewise affine functions on the nearby grid $(0,3h/2,2h,3h,\dots,1)$ that vanish at 0 and 1.  Let  
\[
a_h^+(u,w)=a_h(u,w) = \int_0^1 \frac{\partial u}{\partial x} \frac{\partial w}{\partial x} \, dx,
\]
so that the projectors $r_h$ and $r_h^+$ coincide with the nodal interpolants onto $\mathcal{V}_h$ and $\mathcal{V}_h^+$, respectively~\cite[Remark 3.25(i)]{Ern2004}.  In this setting, the conditions of Corollaries~\ref{corollary1} and~\ref{corollary2} hold with $\eta=\infty$, $\gamma=1$, $r=2$, and $\ell(h) \equiv 1$, leading to the estimates
\begin{align*}
\| r_h^+ u - r_h u \|_{0,2} &\le C h^{5/2} |u|_{r,\infty}, \\
\| r_h^+ u - r_h u \|_{1,2} &\le C h^{3/2} |u|_{r,\infty}
\end{align*}
for $u \in W^{2,\infty}(0,1) \cap H^1_0(0,1)$.

However, consider the function
\[
u(x) = x^{2-1/p} - x
\] 
with $2 < p < \infty$, so that $u \in W^{2,p-\varepsilon}(0,1) \cap H^1_0(0,1)$ for any $\varepsilon>0$.  Then a direct calculation renders that
\begin{align*}
\| r_h^+ u - r_h u \|_{0,2} &\ge C h^{5/2-1/p}, \\
\| r_h^+ u - r_h u \|_{1,2} &\ge C h^{3/2-1/p},
\end{align*}
which are of a lower order than the rates $h^{5/2}$ and $h^{3/2}$, respectively, obtainable for a function in $W^{2,\infty}(0,1) \cap H^1_0(0,1)$.  In fact, by letting $p \rightarrow 2$, these rates can be made arbitrarily close to the quadratic and linear rates that hold in the $L^2$- and $H^1$-norms, respectively, on a pair of unrelated meshes. 

\section{Numerical Examples} \label{sec:numerical}

In this section, we numerically illustrate the superconvergent estimates of Corollaries~\ref{corollary1} and~\ref{corollary2} on test cases in one and two dimensions.

\begin{table}
\centering
\caption{$L^2$-supercloseness of $L^2$-projections onto piecewise affine ($r=2$) and piecewise quadratic ($r=3$) finite element spaces over nearby meshes ($\gamma=1$) in one dimension.}
\label{tab:L2proj1d}
\newcolumntype{R}{>{\raggedleft\arraybackslash}X}%
\begin{tabularx}{295pt}{r|cc|cc}
\hline\noalign{\smallskip}
 &  \multicolumn{2}{c|}{Affine ($r=2$)} & \multicolumn{2}{c}{Quadratic ($r=3$)} \\
$h_0/h$ & \multicolumn{1}{c}{$\|r_h^+ u - r_h u\|_{0,2}$} & \multicolumn{1}{c|}{Order} & \multicolumn{1}{c}{$\|r_h^+ u - r_h u\|_{0,2}$} & \multicolumn{1}{c}{Order} \\
\noalign{\smallskip}\hline\noalign{\smallskip}
  1  &  3.2150e-03  &  \multicolumn{1}{c|}{-}  &  1.2843e-04  &  \multicolumn{1}{c}{-} \\ 
  2  &  5.6505e-04  &      2.5084  &  1.0676e-05  &      3.5886 \\ 
  4  &  9.9837e-05  &      2.5007  &  9.1277e-07  &      3.5480 \\ 
  8  &  1.7645e-05  &      2.5003  &  7.9301e-08  &      3.5248 \\ 
 16  &  3.1189e-06  &      2.5002  &  6.9484e-09  &      3.5126 \\ 
 32  &  5.5132e-07  &      2.5001  &  6.1146e-10  &      3.5063 \\ 
\end{tabularx}
\end{table}

\begin{table}
\centering
\caption{$H^1$-supercloseness of elliptic projections onto piecewise affine ($r=2$) and piecewise quadratic ($r=3$) finite element spaces over nearby meshes ($\gamma=1$) in one dimension.}
\label{tab:ellipticproj1dH1}
\newcolumntype{R}{>{\raggedleft\arraybackslash}X}%
\begin{tabularx}{295pt}{r|cc|cc}
\hline\noalign{\smallskip}
 &  \multicolumn{2}{c|}{Affine ($r=2$)} & \multicolumn{2}{c}{Quadratic ($r=3$)} \\
$h_0/h$ & \multicolumn{1}{c}{$\|r_h^+ u - r_h u\|_{1,2}$} & \multicolumn{1}{c|}{Order} & \multicolumn{1}{c}{$\|r_h^+ u - r_h u\|_{1,2}$} & \multicolumn{1}{c}{Order} \\
\noalign{\smallskip}\hline\noalign{\smallskip}
  1  &  1.4451e-01  &  \multicolumn{1}{c|}{-}  &  7.4390e-03  &  \multicolumn{1}{c}{-} \\ 
  2  &  5.1203e-02  &      1.4968  &  1.2835e-03  &      2.5351 \\ 
  4  &  1.8081e-02  &      1.5017  &  2.2408e-04  &      2.5180 \\ 
  8  &  6.3851e-03  &      1.5017  &  3.9364e-05  &      2.5090 \\ 
 16  &  2.2558e-03  &      1.5011  &  6.9369e-06  &      2.5045 \\ 
 32  &  7.9723e-04  &      1.5006  &  1.2243e-06  &      2.5023 \\ 
\end{tabularx}
\end{table}

\begin{table}
\centering
\caption{$L^2$-supercloseness of elliptic projections onto piecewise affine ($r=2$) and piecewise quadratic ($r=3$) finite element spaces over nearby meshes ($\gamma=1$) in one dimension.}
\label{tab:ellipticproj1dL2}
\newcolumntype{R}{>{\raggedleft\arraybackslash}X}%
\begin{tabularx}{295pt}{r|cc|cc}
\hline\noalign{\smallskip}
 &  \multicolumn{2}{c|}{Affine ($r=2$)} & \multicolumn{2}{c}{Quadratic ($r=3$)} \\
$h_0/h$ & \multicolumn{1}{c}{$\|r_h^+ u - r_h u\|_{0,2}$} & \multicolumn{1}{c|}{Order} & \multicolumn{1}{c}{$\|r_h^+ u - r_h u\|_{0,2}$} & \multicolumn{1}{c}{Order} \\
\noalign{\smallskip}\hline\noalign{\smallskip}
  1  &  3.4546e-03  &  \multicolumn{1}{c|}{-}  &  1.7770e-04  &  \multicolumn{1}{c}{-} \\ 
  2  &  6.1937e-04  &      2.4796  &  1.5493e-05  &      3.5198 \\ 
  4  &  1.1019e-04  &      2.4908  &  1.3576e-06  &      3.5124 \\ 
  8  &  1.9537e-05  &      2.4957  &  1.1943e-07  &      3.5069 \\ 
 16  &  3.4587e-06  &      2.4979  &  1.0530e-08  &      3.5036 \\ 
 32  &  6.1186e-07  &      2.4990  &  9.2955e-10  &      3.5018 \\ 
\end{tabularx}
\end{table}

\paragraph{One dimension.}%
Consider the case in which $\mathcal{V}_h$ is the space of piecewise polynomial functions of degree at most $r-1$ on a grid $(0,h,2h,3h,\dots,1)$ of the unit interval in one dimension that vanish at 0 and 1. Let $\mathcal{V}_h^+$ be the space of piecewise polynomial functions of the same degree that vanish at 0 and 1, on the same grid but with the node nearest to $x=1/4$ perturbed by $h/4$ in the positive direction.  In this scenario, assumption~(\ref{assumption2b}) is satisfied with $\gamma=1$.  Let $u(x) = \sin(\pi x)$ and let 
\[
a_h^+(u,w)=a_h(u,w) = \int_0^1 uw \, dx,
\]
so that $r_h$ and $r_h^+$ are the $L^2$-projectors onto $\mathcal{V}_h$ and $\mathcal{V}_h^+$, respectively.  

Table~\ref{tab:L2proj1d} shows the $L^2$-norm of the difference $r_h^+ u - r_h u$ for several values of $h$, beginning with $h = 1/8 =: h_0$.  The table illustrates the predictions of Corollary~\ref{corollary1}, namely
\begin{equation*}
\|r_h^+ u - r_h u\|_{0,2} \le 
\begin{cases}
C h^{5/2} |u|_{2,\infty} &\mbox{if } r=2, \\
C h^{7/2} |u|_{3,\infty} &\mbox{if } r = 3.
\end{cases}
\end{equation*}

Next, consider the same setup as above, but with
\[
a_h^+(u,w)=a_h(u,w) = \int_0^1 \frac{\partial u}{\partial x} \frac{\partial w}{\partial x} \, dx,
\]
so that $r_h$ and $r_h^+$ are the standard elliptic projectors onto $\mathcal{V}_h$ and $\mathcal{V}_h^+$, respectively.  Table~\ref{tab:ellipticproj1dH1} shows the $H^1$ norm of the difference $r_h^+ u - r_h u$ for the sequence of grids described above.  The table illustrates the predictions of Corollary~\ref{corollary1}, namely
\begin{equation*}
\|r_h^+ u - r_h u\|_{1,2} \le 
\begin{cases}
C h^{3/2} \log(h^{-1}) |u|_{2,\infty} &\mbox{if } r=2, \\
C h^{5/2} |u|_{3,\infty} &\mbox{if } r = 3.
\end{cases}
\end{equation*}
Table~\ref{tab:ellipticproj1dL2} shows the $L^2$-norm of the difference $r_h^+ u - r_h u$ for the same sequence of grids.  The table illustrates the predictions of Corollary~\ref{corollary2}, namely
\begin{equation*}
\|r_h^+ u - r_h u\|_{0,2} \le 
\begin{cases}
C h^{5/2} \log(h^{-1}) |u|_{2,\infty} &\mbox{if } r=2, \\
C h^{7/2} |u|_{3,\infty} &\mbox{if } r = 3.
\end{cases}
\end{equation*}
Note that we have not attempted to detect the presence of the factor $\log(h^{-1})$ in these numerical experiments.

\begin{figure}[t]
  \centering
  \subfigure[]{
    \centering
    \includegraphics[trim = 1.2in 0.4in 1.1in 0.3in, clip=true, width=0.3\textwidth]{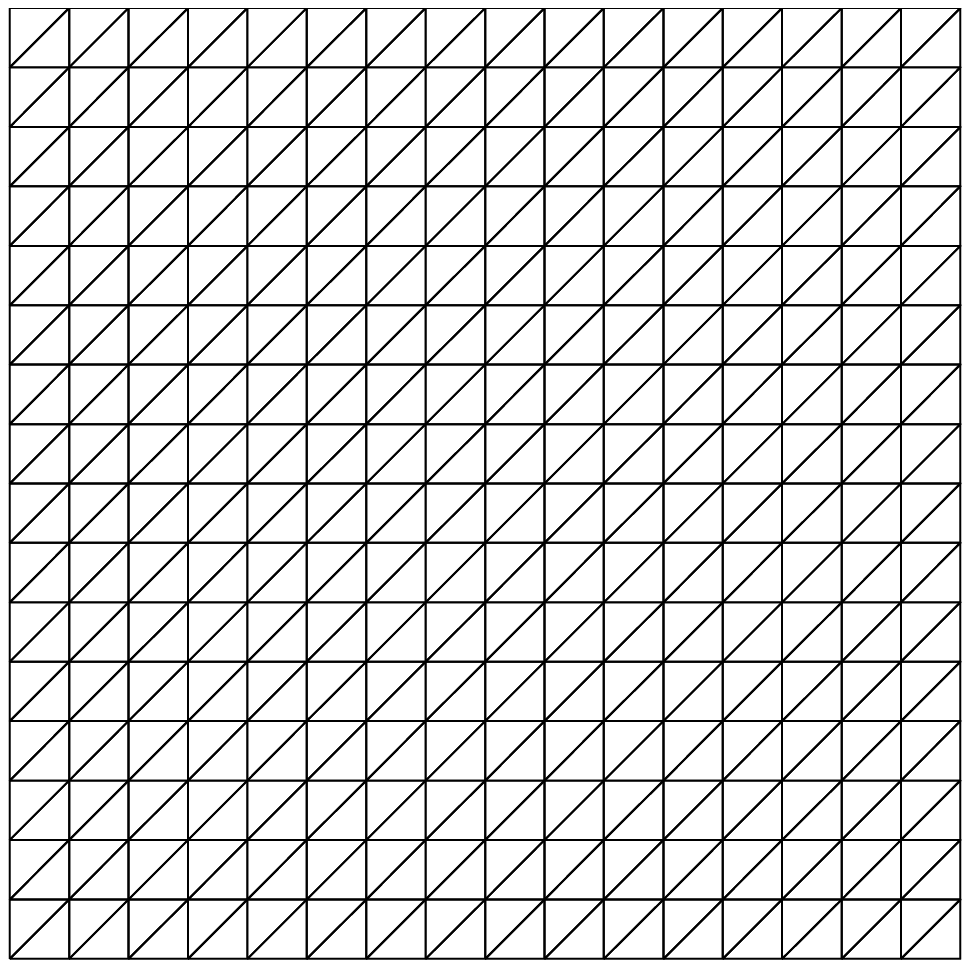}
    \label{fig:squaremesh0}
  }
  \subfigure[]{
    \centering
    \includegraphics[trim = 1.2in 0.4in 1.1in 0.3in, clip=true, width=0.3\textwidth]{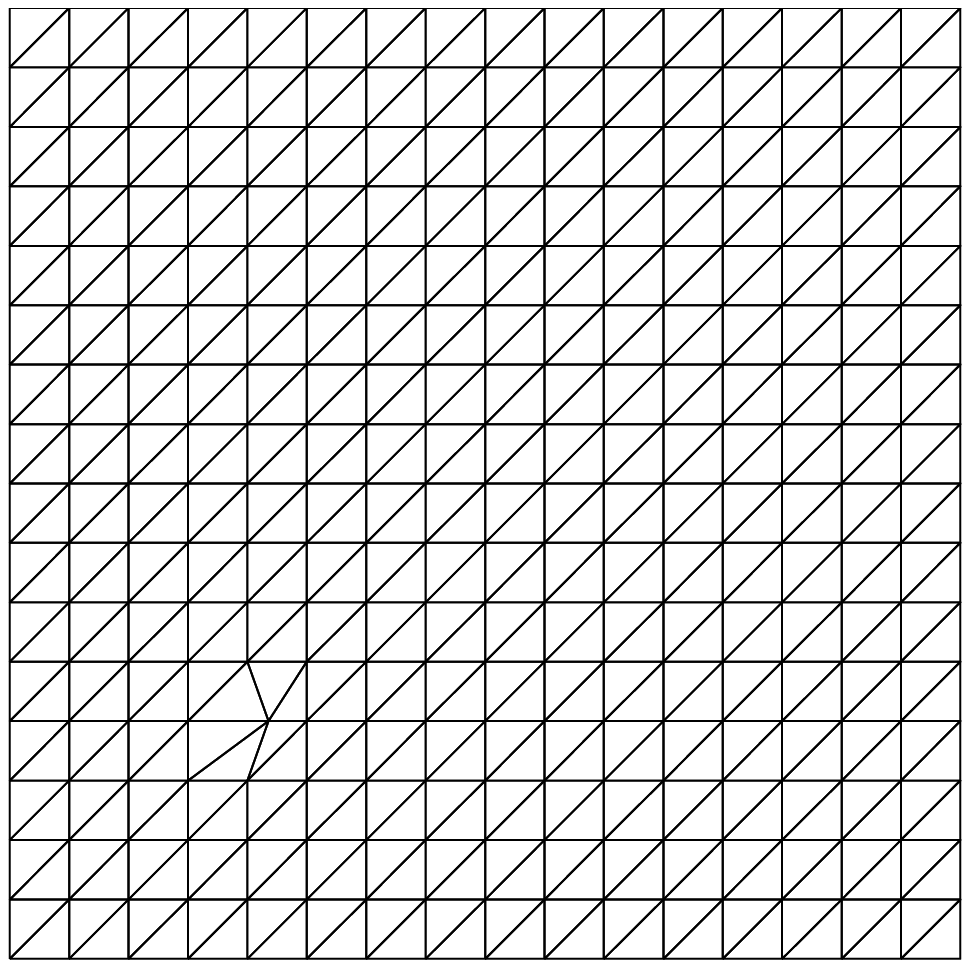}
    \label{fig:squaremesh1}
  }      
  \subfigure[]{
    \centering
    \includegraphics[trim = 1.2in 0.4in 1.1in 0.3in, clip=true, width=0.3\textwidth]{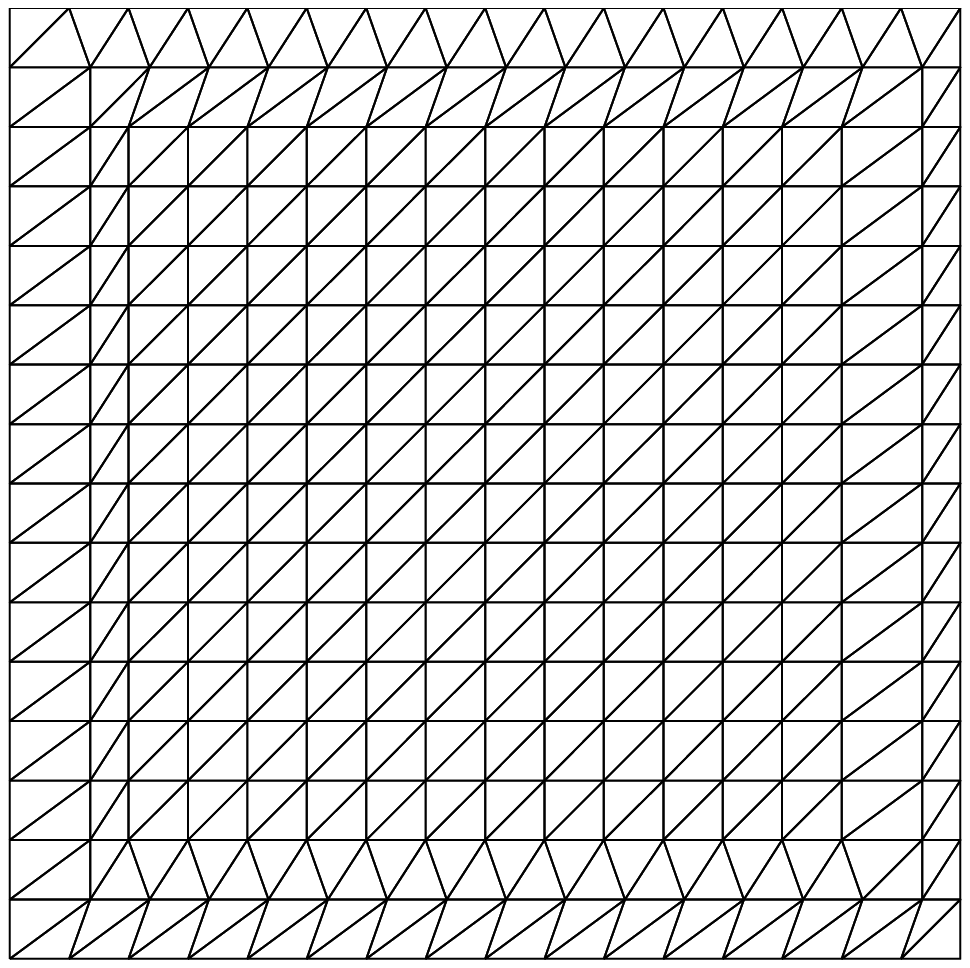}
    \label{fig:squaremesh2}
  }      
  \caption{(a) Mesh of the unit square consisting of equally sized isosceles right triangles. (b) Identical mesh, but with the node at $(x,y)=(1/4,1/4)$ perturbed by $h/4$ in the positive $x$ direction. (c) Identical mesh, but with all nodes having distance $h/\sqrt{2}$ from the boundary perturbed by $h/4$ in the positive $x$ direction.}
  \label{fig:squaremesh}
\end{figure}

\begin{table}[t]
\centering
\caption{$L^2$-supercloseness of $L^2$-projections onto piecewise affine ($r=2$) finite element spaces over nearby meshes ($\gamma=2$; see Figs.~\ref{fig:squaremesh}(a) and~\ref{fig:squaremesh}(b)) in two dimensions.}
\label{tab:L2proj2d}
\newcolumntype{R}{>{\raggedleft\arraybackslash}X}%
\begin{tabularx}{165pt}{r|cc}
\hline\noalign{\smallskip}
 &  \multicolumn{2}{c}{Affine ($r=2$)} \\
$h_0/h$ & \multicolumn{1}{c}{$\|r_h^+ u - r_h u\|_{0,2}$} & \multicolumn{1}{c}{Order} \\
\noalign{\smallskip}\hline\noalign{\smallskip}
  1  &  6.3533e-03  &  \multicolumn{1}{c}{-} \\ 
  2  &  7.5614e-04  &      3.0708 \\ 
  4  &  8.8718e-05  &      3.0914 \\ 
  8  &  1.1020e-05  &      3.0091 \\ 
 16  &  1.3781e-06  &      2.9993 \\ 
\end{tabularx}
\end{table}

\begin{table}[t]
\centering
\caption{$H^1$- and $L^2$-supercloseness of elliptic projections onto piecewise affine ($r=2$) finite element spaces over nearby meshes ($\gamma=2$; see Figs.~\ref{fig:squaremesh}(a) and~\ref{fig:squaremesh}(b)) in two dimensions.}
\vspace{0.1in}
\label{tab:ellipticproj2d}
\newcolumntype{R}{>{\raggedleft\arraybackslash}X}%
\begin{tabularx}{295pt}{r|cc|cc}
\hline\noalign{\smallskip}
 &  \multicolumn{4}{c}{Affine ($r=2$)} \\
$h_0/h$ & \multicolumn{1}{c}{$\|r_h^+ u - r_h u\|_{1,2}$} & \multicolumn{1}{c|}{Order} & \multicolumn{1}{c}{$\|r_h^+ u - r_h u\|_{0,2}$} & \multicolumn{1}{c}{Order} \\
\noalign{\smallskip}\hline\noalign{\smallskip}
  1  &  2.1441e-01  &  \multicolumn{1}{c|}{-}  &  6.6386e-03  &  \multicolumn{1}{c}{-} \\ 
  2  &  4.7374e-02  &      2.1782  &  7.8678e-04  &      3.0768 \\ 
  4  &  1.1359e-02  &      2.0603  &  9.6370e-05  &      3.0293 \\ 
  8  &  2.8114e-03  &      2.0144  &  1.2033e-05  &      3.0016 \\ 
 16  &  7.0176e-04  &      2.0023  &  1.5106e-06  &      2.9937 \\ 
\end{tabularx}
\end{table}

\paragraph{Two dimensions.}%
Consider now the case in which $\mathcal{V}_h \subset H^1_0((0,1) \times (0,1))$ is the space of piecewise affine functions on a mesh of the unit square in two dimensions consisting of equally sized isosceles right triangles, as in Fig.~\ref{fig:squaremesh}(a). Let $\mathcal{V}_h^+ \subset H^1_0((0,1) \times (0,1))$ be the space of piecewise affine functions on the same mesh, but with the node nearest to $(x,y)=(1/4,1/4)$ perturbed by $h/4$ in the positive $x$ direction, as in Fig.~\ref{fig:squaremesh}(b).  In this scenario, assumption~(\ref{assumption2b}) is satisfied with $\gamma=2$.  Let $u(x) = \sin(\pi x)\sin(\pi y)$ and let 
\[
a_h^+(u,w)=a_h(u,w) = \int_0^1\int_0^1 uw \, dx dy,
\]
so that $r_h$ and $r_h^+$ are the $L^2$-projectors onto $\mathcal{V}_h$ and $\mathcal{V}_h^+$, respectively.

Table~\ref{tab:L2proj2d} shows the $L^2$-norm of the difference $r_h^+ u - r_h u$ for several values of $h$, beginning with $h = \sqrt{2}/4 =: h_0$.  The table illustrates the predictions of Corollary~\ref{corollary1}, namely
\begin{equation} \label{est1}
\|r_h^+ u - r_h u\|_{0,2} \le C h^3 |u|_{2,\infty}.
\end{equation}

Next, consider the same setup as above, but with
\[
a_h^+(u,w)=a_h(u,w) = \int_0^1\int_0^1 \left( \frac{\partial u}{\partial x} \frac{\partial w}{\partial x} + \frac{\partial u}{\partial y} \frac{\partial w}{\partial y}  \right) \, dx dy,
\]
so that $r_h$ and $r_h^+$ are the elliptic projectors onto $\mathcal{V}_h$ and $\mathcal{V}_h^+$, respectively.  Table~\ref{tab:ellipticproj2d} shows the $H^1$- and $L^2$-norms of the difference $r_h^+ u - r_h u$ for the sequence of meshes described above.  The table illustrates the predictions of Corollaries~\ref{corollary1} and~\ref{corollary2}, namely
\begin{equation} \label{est2}
\|r_h^+ u - r_h u\|_{m,2} \le 
\begin{cases}
C h^2 \log(h^{-1}) |u|_{2,\infty} &\mbox{if } m=0, \\
C h^3 \log(h^{-1}) |u|_{2,\infty} &\mbox{if } m=1.
\end{cases}
\end{equation}
Again, we have not attempted to detect the presence of the factor $\log(h^{-1})$.

\begin{table}[t]
\centering
\caption{$L^2$-supercloseness of $L^2$-projections onto piecewise affine ($r=2$) finite element spaces over nearby meshes ($\gamma=1$; see Figs.~\ref{fig:squaremesh}(a) and~\ref{fig:squaremesh}(c)) in two dimensions.  Relative to Table~\ref{tab:L2proj2d}, a lower order of superconvergence is observed due to the larger fraction of perturbed elements present in the perturbed mesh.}
\label{tab:L2proj2dgamma1}
\newcolumntype{R}{>{\raggedleft\arraybackslash}X}%
\begin{tabularx}{165pt}{r|cc}
\hline\noalign{\smallskip}
 &  \multicolumn{2}{c}{Affine ($r=2$)} \\
$h_0/h$ & \multicolumn{1}{c}{$\|r_h^+ u - r_h u\|_{0,2}$} & \multicolumn{1}{c}{Order} \\
\noalign{\smallskip}\hline\noalign{\smallskip}
  1  &  2.2504e-02  &  \multicolumn{1}{c}{-} \\ 
  2  &  4.8445e-03  &      2.2158 \\ 
  4  &  1.0019e-03  &      2.2736 \\ 
  8  &  1.9159e-04  &      2.3866 \\ 
 16  &  3.5132e-05  &      2.4472 \\ 
 32  &  6.3195e-06  &      2.4749 \\
\end{tabularx}
\end{table}

\begin{table}[t]
\centering
\caption{$H^1$- and $L^2$-supercloseness of elliptic projections onto piecewise affine ($r=2$) finite element spaces over nearby meshes ($\gamma=1$; see Figs.~\ref{fig:squaremesh}(a) and~\ref{fig:squaremesh}(c)) in two dimensions.  Relative to Table~\ref{tab:ellipticproj2d}, lower orders of superconvergence are observed due to the larger fraction of perturbed elements present in the perturbed mesh.}
\vspace{0.1in}
\label{tab:ellipticproj2dgamma1}
\newcolumntype{R}{>{\raggedleft\arraybackslash}X}%
\begin{tabularx}{295pt}{r|cc|cc}
\hline\noalign{\smallskip}
 &  \multicolumn{4}{c}{Affine ($r=2$)} \\
$h_0/h$ & \multicolumn{1}{c}{$\|r_h^+ u - r_h u\|_{1,2}$} & \multicolumn{1}{c|}{Order} & \multicolumn{1}{c}{$\|r_h^+ u - r_h u\|_{0,2}$} & \multicolumn{1}{c}{Order} \\
\noalign{\smallskip}\hline\noalign{\smallskip}
  1  &  5.4318e-01  &  \multicolumn{1}{c|}{-}  &  1.9864e-02  &  \multicolumn{1}{c}{-} \\ 
  2  &  2.8504e-01  &      0.9303  &  4.8794e-03  &      2.0254 \\ 
  4  &  1.2522e-01  &      1.1867  &  1.0528e-03  &      2.2125 \\ 
  8  &  4.8674e-02  &      1.3632  &  1.9842e-04  &      2.4075 \\ 
 16  &  1.7931e-02  &      1.4407  &  3.5671e-05  &      2.4758 \\ 
 32  &  6.4595e-03  &      1.4730  &  6.3290e-06  &      2.4947 \\
\end{tabularx}
\end{table}

\paragraph{More substantial mesh perturbation in two dimensions.}%
Finally, consider the same two-dimensional tests as above, but with the mesh of Fig.~\ref{fig:squaremesh}(b) replaced by a different perturbation of the uniform mesh.  Namely, consider perturbing all nodes whose distance from the boundary of the unit square is equal to $h/\sqrt{2}$ (the length of the shortest edge of each triangle) via a translation by $h/4$ in the positive $x$ direction, as in Fig.~\ref{fig:squaremesh}(c).

In this scenario, assumption~(\ref{assumption2b}) is satisfied with $\gamma=1$, so that the estimates~(\ref{est1}) and~(\ref{est2}) no longer apply.  Their analogues in this case read
\begin{equation*}
\|r_h^+ u - r_h u\|_{0,2} \le C h^{5/2} |u|_{2,\infty}.
\end{equation*}
and
\begin{equation*}
\|r_h^+ u - r_h u\|_{m,2} \le 
\begin{cases}
C h^{3/2} \log(h^{-1}) |u|_{2,\infty} &\mbox{if } m=0, \\
C h^{5/2} \log(h^{-1}) |u|_{2,\infty} &\mbox{if } m=1,
\end{cases}
\end{equation*}
respectively.  Tables~\ref{tab:L2proj2dgamma1}-\ref{tab:ellipticproj2dgamma1} illustrate these predictions.  Again, we have not attempted to detect the presence of the factor $\log(h^{-1})$.

\section{Summary} \label{sec:conclusion}

We have derived estimates for the difference between the orthogonal projections $r_h u$ and $r_h^+ u$ of a smooth function $u$ onto nearby finite element spaces $\mathcal{V}_h$ and $\mathcal{V}_h^+$, respectively, with respect to bilinear forms $a_h, a_h^+ : \mathcal{V} \times \mathcal{V} \rightarrow \mathbb{R}$, respectively, where $\mathcal{V}$ is a closed subspace of $H^s(\Omega)$. When $s \in \{0,1\}$ and $\mathcal{V}_h$ and $\mathcal{V}_h^+$ consist of continuous functions that are elementwise polynomials over shape-regular, quasi-uniform meshes that coincide except on a region of measure $O(h^\gamma)$ for a constant $\gamma \ge 0$, the estimates for $\|r_h^+ u - r_h u\|_{s,2}$ are superconvergent by $O(h^{\gamma/2})$, provided that $u \in W^{s,\infty}(\Omega)$ and $a_h$ and $a_h^+$ are sufficiently close.  In addition, when $s=1$ and a few more mild assumptions (namely~(\ref{cond:smoothing}-\ref{cond:interp_coincide})) are satisfied, an $O(h^{\gamma/2})$-superconvergent estimate for $\|r_h^+ u - r_h u\|_{0,2}$ holds.  Numerical experiments illustrated these estimates and verified the necessity of the regularity assumptions on $u$.

\section{Acknowledgments}
This research was supported by the U.S. Department of Energy, grant number DE-FG02-97ER25308; Department of the Army Research Grant, grant number: W911NF-07- 2-0027; NSF Career Award, grant number: CMMI-0747089; and NSF, grant number CMMI-1301396.

\bibliographystyle{unsrtnat}
\bibliography{projections}{}

\appendix 
\section{Properties of Piecewise Polynomial Finite Element Spaces} \label{sec:appendix_Pk}

In this section, we verify conditions~(\ref{Wseta}-\ref{assumption2b}) for piecewise polynomial finite element spaces on nearby meshes for the cases $s=0$ and $s=1$.

As in Section~\ref{sec:intro}, consider two families of shape-regular, quasi-uniform meshes $\{\mathcal{T}_h\}_{h \le h_0}$ and $\{\mathcal{T}_h^+\}_{h \le h_0}$ of an open, bounded, Lipschitz domain $\Omega \subset \mathbb{R}^d$, $d \ge 1$.  Assume that the two families are parametrized by a scalar $h$ that equals the maximum diameter of an element among all elements of $\mathcal{T}_h$ and $\mathcal{T}_h^+$ for every $h \le h_0$.  Let $\mathcal{V}_h$ and $\mathcal{V}_h^+$ be finite element spaces consisting of continuous functions that are elementwise polynomials of degree at most $r-1$ over $\mathcal{T}_h$ and $\mathcal{T}_h^+$, respectively, where $r>1$ is an integer.

In this setting, condition~(\ref{Wseta}) is automatic for any $\eta \in [2,\infty]$, $s \in \{0,1\}$.  Condition~(\ref{inverse}) is trivial for $s=0$ and is satisfied for $s=1$ and any $\eta \in [2,\infty]$~\cite{Ern2004}.  

Condition~(\ref{assumption2b}) holds for any $\eta \in [2,\infty]$ when $\mathcal{T}_h$ and $\mathcal{T}_h^+$ coincide except on a region of measure $O(h^\gamma)$.  To prove this, let $\{N_a\}_{a=1}^A \subset \mathcal{V}_h$ and $\{N_a^+\}_{a=1}^{A^+} \subset \mathcal{V}_h^+$ be the standard Lagrange shape functions that form bases for $\mathcal{V}_h$ and $\mathcal{V}_h^+$, respectively.  Our assumptions on $\mathcal{T}_h$ and $\mathcal{T}_h^+$ imply the existence of an integer $I$ such that $N_a = N_a^+$ for every $1 \le a \le I$ and such that
\begin{equation}
\left| \left( \displaystyle\bigcup_{a=I+1}^A \mathrm{supp}(N_a) \right) \cup \left( \displaystyle\bigcup_{a=I+1}^{A^+} \mathrm{supp}(N_a^+) \right) \right| \le C h^\gamma
\end{equation}
for every $h \le h_0$.  

Define $\pi_h : \mathcal{V}_h^+ + \mathcal{V}_h \rightarrow \mathcal{V}_h^+ \cap \mathcal{V}_h$ as follows: For any
\begin{equation} \label{wh_expansion}
w_h = \sum_{a=1}^I c_a N_a + \sum_{a=I+1}^A c_a N_a + \sum_{a=I+1}^{A^+} c_a^+ N_a^+
\end{equation}
belonging to $\mathcal{V}_h^+ + \mathcal{V}_h$, set
\begin{equation} \label{pih_def}
\pi_h w_h := \sum_{a=1}^I c_a N_a.
\end{equation}

Clearly,
\[
|\mathrm{supp}(\pi_h w_h - w_h)| \le C h^\gamma
\]
for every $w_h \in \mathcal{V}_h^+ + \mathcal{V}_h$ and every $h \le h_0$.  To prove that
\begin{equation} \label{pih_stability}
\| \pi_h w_h \|_{0,\eta} \le C \|w_h\|_{0,\eta}
\end{equation}
for every $w_h \in \mathcal{V}_h^+ + \mathcal{V}_h$ and every $h \le h_0$, there are two cases to consider: $\eta = \infty$ and $2 \le \eta < \infty$.

For $\eta = \infty$, it is enough to note that for each of the two finite element spaces, every shape function is bounded uniformly in $h$ in the maximum norm, the number of shape functions whose support intersects any given element is bounded uniformly in $h$, and the coefficients $c_a$, $1 \le a \le I$, in the expansion~(\ref{wh_expansion}) of $w_h$ are bounded by $\|w_h\|_{0,\infty}$.  Indeed, the standard degrees of freedom $\sigma_a$, $1 \le a \le I$, for the Lagrange shape functions $N_a (=N_a^+)$, $1 \le a \le I$, satisfy
\[
\sigma_a(N_b) = \delta_{ab}, \;\;\; 1 \le b \le A
\]
and
\[
\sigma_a(N_b^+) = \delta_{ab}, \;\;\; 1 \le b \le A^+,
\]
where $\delta_{ab}$ denotes the Kronecker delta.
Hence, for any $1 \le a \le I$, 
\[
|c_a| = |\sigma_a(w_h)| \le \|w_h\|_{0,\infty}.
\]

For $2 \le \eta < \infty$, the proof of~(\ref{pih_stability}) relies on the following lemma.

\begin{lemma} \label{lemma:norms}
Let $\{\mathcal{T}_h\}_{h \le h_0}$ be a shape-regular, quasi-uniform family of meshes of an open, bounded, Lipschitz domain $\Omega \subset \mathbb{R}^d$, $d \ge 1$, with $h$ denoting the maximum diameter of an element $K \in \mathcal{T}_h$.  Let $r>1$ be an integer.  For any $K \in \mathcal{T}_h$, let $\theta_1,\theta_2,\dots,\theta_{n_{sh}}$ denote the local shape functions for the Lagrange finite element of degree at most $r-1$ on $K$.  Then for any $2 \le \eta < \infty$, there exist $C_1,C_2 > 0$ independent of $h$ such that for every $h \le h_0$, every $K \in \mathcal{T}_h$, and every $v = \sum_{i=1}^{n_{sh}} d_i \theta_i$,
\[
C_1 h^d \sum_{i=1}^{n_{sh}} |d_i|^\eta \le \|v\|_{0,\eta,K}^\eta \le C_2 h^d \sum_{i=1}^{n_{sh}} |d_i|^\eta.
\]
\end{lemma}
\proof A proof of this fact when $\eta=2$ is given in~\cite[Lemma 9.7]{Ern2004}.  The case $2 < \eta < \infty$ is a trivial modification thereof.  
\qed

Now let $w_h$ and $\pi_h w_h$ be as in~(\ref{wh_expansion}) and~(\ref{pih_def}), respectively.  Note that the support of $\pi_h w_h$ is contained within the region $Q_h \subseteq \Omega$ over which $\mathcal{T}_h$ and $\mathcal{T}_h^+$ coincide.  On any $K \in \mathcal{T}_h$ with $K \subseteq Q_h$, we can write 
\[
\left.w_h\right|_K = \sum_{i=1}^{n_{sh}} d_i \theta_i
\]
and
\[
\left.\pi_h w_h\right|_K = \sum_{i=1}^{n_{sh}} \bar{d}_i \theta_i,
\]
with scalars $d_i \in \mathbb{R}$ and $\bar{d}_i \in \{0,d_i\}$ for every $i$.  By Lemma~\ref{lemma:norms},
\begin{align*}
\| \pi_h w_h \|_{0,\eta,K}^\eta 
&\le C_2 h^d \sum_{i=1}^{n_{sh}} |\bar{d}_i|^\eta \\
&\le C_2 h^d \sum_{i=1}^{n_{sh}} |d_i|^\eta \\
&\le C_2 C_1^{-1}  \| w_h \|_{0,\eta,K}^\eta
\end{align*}
on every such $K$.  Summing over all $K \in \mathcal{T}_h$ with $K \subseteq Q_h$ proves~(\ref{pih_stability}) for $2 \le \eta < \infty$.

\section{Estimates for the $L^2$-Projection and Elliptic Projections} \label{sec:appendix_Linf}
Two exemplary cases in which estimates of the form~(\ref{ep0eta}-\ref{epmeta}) are known to hold are the following.  Suppose that $\mathcal{V} = H^s(\Omega) \cap H^1_0(\Omega)$ and $\mathcal{V}_h$ is the space of continuous functions in $\mathcal{V}$ that are elementwise polynomials of degree at most $r-1$ on a shape-regular, quasi-uniform family of meshes $\{\mathcal{T}_h\}_{h \le h_0}$ whose maximum element diameter is $h$.  Then:
\begin{enumerate}[label=(\roman*)]
\item If $s=0$, $d \in \{1,2\}$, and 
\[a_h(u,w) = \int_\Omega uw \, dx
\]
so that $r_h$ is the $L^2$-projector onto $\mathcal{V}_h$, then~(\ref{ep0eta}) holds with $\ell(h) \equiv 1$ for any $\eta \in [2,\infty]$~\cite{Crouzeix1987}.  Note that the estimate~(\ref{epmeta}) is vacuous in this case, since $s=0$.
\item If $s=1$, $d \in \{2,3\}$, and 
\[
a_h(u,w) = \int_{\Omega} \left(\sum_{i,j=1}^d a_{ij}(x) \frac{\partial u}{\partial x_i} \frac{\partial w}{\partial x_j} + \sum_{j=1}^d b_j(x) \frac{\partial u}{\partial x_j} w + b_0(x) u w\right) dx
\] 
with $h$-independent coefficients $a_{ij}$, $i,j=1,2,\dots,d$ and $b_j$, $j=0,1,\dots,d$, then~(\ref{ep0eta}-\ref{epmeta}) hold~\cite{Ern2004} with $\ell(h) \equiv 1$ for any $2 \le \eta < \infty$ (if $r=2$) and any $\eta \in [2,\infty]$ (if $r>2$),
provided that
\begin{itemize}
\item The coefficients satisfy $b_j \in L^\infty(\Omega)$, $j=0,1,\dots,d$, and $a_{ij} \in L^\infty(\Omega) \cap W^{1,p}(\Omega)$, $i,j=1,2,\dots,d$, with $p>2$ if $d=2$ and $p \ge 12/15$ if $d=3$.
\item The coefficients $a_{ij}$ are coercive pointwise, i.e. there exists $c>0$ independent of $x$ such that
\begin{equation} \label{pointwise_coercive}
\sum_{i,j=1}^d a_{ij}(x) \xi_i \xi_j \ge c |\xi|^2
\end{equation}
for every $0 \ne \xi \in \mathbb{R}^d$ and a.e. $x \in \Omega$.
\item There exists $C>0$, $q_0 > d$ such that the continuous Dirichlet problem
\[
a_h(u,w) = \int_\Omega f w \, dx \;\;\; \forall w \in \mathcal{V}
\]
has a unique solution satisfying
\begin{equation} \label{dirichlet_bound}
\|u\|_{2,q} \le C \|f\|_{0,q}
\end{equation}
for every $f \in L^p(\Omega)$ and every $1 < q < q_0$.
\end{itemize}
Under the same conditions as above but with $r=2$ and $\eta=\infty$, the estimates~(\ref{ep0eta}-\ref{epmeta}) hold with $\ell(h) = \log(h^{-1})$ in dimension $d=2$~\cite{Ern2004}.
\end{enumerate}

\end{document}